\documentclass[12pt]{article}
\usepackage{mathrsfs}
\usepackage{amsfonts}
\usepackage{bbm}
\usepackage{amssymb}
\usepackage{amsmath,amssymb}
\usepackage{latexcad}

\newcommand{\eset}{\emptyset}

\newcommand{\G}{\Gamma}

\newcommand{\Llra}{\Longleftrightarrow}

\newcommand{\mf}{\mathfrak}
\newcommand{\nin}{\noindent}

\newcommand{\om}{\omega}
\newcommand{\Om}{\Omega}
\newcommand{\ol}{\overline}
\newcommand{\os}{\overset}

\newcommand{\seq}{\subseteq}

\newcommand{\vs}{\vspace*}
\newcommand{\vi}{\varphi}

\def\Lra{\Longrightarrow}
\def \nin {\noindent}

\def \Lemma #1 {\vs{3mm}\nin {\bf Lemma #1} \it}
\def \Prop #1 {\vs{3mm}\nin {\bf Proposition #1} \it}
\def \Th #1 {\vs{3mm}\nin {\bf Theorem #1} \it}
\def \Cor #1 {\vs{3mm}\nin {\bf Corollary #1} \it}
\def \Ex #1 {\vs{3mm}\nin {\bf Example #1} \it}
\def \Proof {\vs{3mm}\nin {\bf Proof. }}
\def \part #1 {\hfil\break\hglue 12pt {\rm (#1)~}}

\def \qed {~\vrule height6pt width 6pt depth 0pt}
\def\fs{\footnotesize}

\title{
\bf\LARGE On graphs related to co-maximal ideals of a commutative ring\thanks {This research is
supported by the National Natural Science Foundation of China (Grant
No.10671122). The second author is also partly supported by STCSM No. 09XD1402500. } }
\author{{Tongsuo Wu\thanks{tswu@sjtu.edu.cn}, Meng Ye\thanks{yebeibei\_1982cn@yahoo.com.cn}
 }\\
 {\small Department of Mathematics, Shanghai Jiaotong University}\\
{\small Shanghai 200240, P. R. China}\\
Dancheng Lu\thanks{
ludancheng@suda.edu.cn}\\
 {\fs Department of Mathematics, Suzhou University,  Suzhou 215006, P.R. China}\\
Houyi Yu\thanks{yhy178@163.com}\\{\small Department of Mathematics, Shanghai Jiaotong University}\\
{\small Shanghai 200240, P. R. China}\\
}

\date{}
\begin{document}
\baselineskip=16pt \maketitle

\begin{center}
\begin{minipage}{12cm}

 \vs{3mm}\nin{\small\bf Abstract.} {\fs This paper studies the co-maximal graph $\Om(R)$, the induced subgraph $\G(R)$ of $\Om(R)$ whose vertex set is $R\setminus (U(R)\cup J(R))$ and a retract $\G_r(R)$ of $\G(R)$, where $R$ is a commutative ring.  We show that the core of $\G(R)$ is a union of triangles and rectangles, while a vertex in $\G(R)$ is either an end vertex or a vertex in the core. For a non-local ring $R$, we prove that both the chromatic number and clique number of $\G(R)$ are identical with the number of maximal ideals of $R$.  A graph $\G_r(R)$ is also introduced on the vertex set $\{Rx|\,x\in R\setminus (U(R)\cup J(R))\}$, and graph properties of $\G_r(R)$ are studied. }

\vs{3mm}\nin {\small Key Words:} {\small Co-maximal graph; Split graph; Core; Chromatic number; Retract of a graph.
}

\end{minipage}
\end{center}

\vs{4mm}\nin{\bf 1. Introduction}

\vs{3mm}\nin  In 1988, Beck \cite{B} introduced the concept of zero-divisor graph for a commutative ring. Since then a lot of work was done in this area of research. Several other graph structures were also defined on rings and semigroups. In 1995, Sharma and Bhatwadekar \cite{SB} introduced a graph $\Om(R) $ on a commutative
ring $R$, whose vertices are elements of R where two distinct vertices
$x$ and $y$ are adjacent if and only if $Rx + Ry = R$. Recently, Maimani et.al. in \cite{MSSY}
named this graph $\Om(R)$ as the {\it co-maximal graph} of $R$ and they noticed that the  subgraph $\G(R)$ induced on the subset $R\setminus (U(R)\cup J(R))$ is the key to the co-maximal graph . Many interesting results about the subgraph were obtained in \cite{MSSY} and Wang \cite{wang}, and their work show that the properties of the graph $\G(R)$ are quite similar to that of the modified zero-divisor graph  by Anderson and Livingston \cite{AL}. For example, both graphs are simple, connected and with diameter less than or equal to three, and each has girth less than or equal to four if they contain a cycle. Because of this reason, in this paper we use $\G(R)$ to denote the graph $\G_2(R)\setminus J(R)$ of \cite{MSSY}. We discover more properties shared by both zero-divisor graph and the subgraph $\G(R)$ of $\Om(R)$. In particular,  It is shown that the core of $\G(R)$ is a union of triangles and rectangles, while a vertex in $\G(R)$ is either an end vertex or a vertex in the core. For any non-local ring $R$, it is shown that the chromatic number of the graph $\G(R)$ is identical with the number of maximal ideals of $R$. In Section 4, we  introduce a new graph $\G_r(R)$ on the vertex set $$\{Rx|\,x\in R\setminus (U(R)\cup J(R))\}.$$ This graph is in fact a retract of the graph $\G(R)$ and thus simpler than the graph $\G(R)$ in general, but we will show that they  share many common properties and invariants.

Jinnah and Mathew in \cite{JM} studied the problem of when a co-maximal graph $\Om(R)$ is a split graph, and they determined all rings $R$ with the property. In Section 2, we give an alternative proof to their Theorem 2.3.
In the co-maximal graph $\Om(R)$, each  unit $u$ of $R$ is adjacent to all vertices of the graph while an element of $J(R)$ only connects to units of $R$.
Temporally, we say $u$ is in the center of the graph  $\G(R)$. Related to the co-maximal relation, there is the concept of {\it rings with stable range one}. Recall that a ring $R$ (which needs not
 be commutative) has one in its stable range, if for any $x,y$ with $Rx+Ry=R$, there is an element $t$ such that $x+ty$ is invertible. For example, the following classes of rings have one in their
stable range: zero-dimensional commutative rings, von Neumann unit-regular rings, semilocal rings. The concept co-maximal graph gives an interesting graph interpretation of such rings. In fact,
{\it  a commutative ring $R$ has one in its stable range if and only if for any pair of adjacent vertices $x,y$ in the co-maximal graph $\Om(R)$, the additive coset $x+Ry$ (and $y+Rx$) has at
 least one element in the center of the graph $\Om(R)$}.

 Throughout this paper, all rings are assumed to be commutative with identity. For a ring $R$,  let $U(R)$ be the set of invertible elements of $R$ and $J(R)$ the Jacobson radical of $R$.  Recall that a graph is called {\it complete} ({\it discrete}, respectively) if every pair of vertices are adjacent (respectively, no pair of vertices are adjacent). We denote a complete graph by $K$,
a complete (discrete, resp.) graph with $n$ vertices  by $K_n$ (resp., $D_n$).  A subset $K$ of the vertex set of $G$ is called a {\it clique} if any two
distinct vertices of $K$ are adjacent; the {\it clique number} $\om(G)$ of $G$ is the least upper bound of the size of the
cliques. Similarly, we denote by $K_{m,n}$ the complete bipartite graph with two partitions of sizes $m,n$ respectively. Recall that a simple graph $G$ is called a {\it refinement}
of a simple graph $H$ if $V(G)=V(H)$ and $a-b$ in $H$ implies $a-b$ in $G$ for all distinct vertices of $G$, where $a-b$ means that $a\not=b$ and $a$ is adjacent to $b$. Recall that a {\it cycle} in a graph is a path $v_1-v_2-\cdots-v_n$ together with an additional edge $v_n-v_1$ ($n\ge 3$).
For a simple graph $G$ and a nonempty subset $S$ of $V(G)$, there is {\it the subgraph induced on $S$ }: the vertex set is $S$ and the edge set is $$\{x-y\,|\,x\not=y\in S,\, \text{and there is an edge $x-y$ in the graph $G$ }\}.$$
A discrete induced  subgraph of a graph $G$ is also called an {\it independent subset} of $G$.

\vs{4mm}\nin{\bf 2. Rings $R$ whose co-maximal graph is a split graph}

\vs{3mm}\nin Throughout this section, assume that $G$ is a {\it split graph}, i.e., $G$ is simple and connected with $V(G)=K\cup D$, where $K\cap D=\eset$ and the induced subgraph on $K$ (respectively, on $D$) is a complete (discrete, respectively) graph. For the split graph $G$, we always assume that $D$ is a maximal such independent subset. Under the assumption, a complete graph $K_n$ is a split graph with $|K|=n-1,\,|D|=1$.

\vs{3mm}\nin{\bf Lemma 2.1.}  {\it For a commutative ring $R$, let $G$ be the co-maximal graph $\Om(R)$. If $G$ is a split graph with $V(G)=K\cup D$, then }

(1) {\it For any proper ideal $I$ of $R$, $|I\cap K|\le 1$.}

(2) {\it $\mf m\cap\mf n\cap K=\{0\}$ holds for distinct maximal ideals $\mf m,\mf n$ of $R$.}

(3) {\it If $R$ is isomorphic to neither $\mathbb Z_2$ nor $\mathbb Z_2\times \mathbb Z_2$, then $|K|\ge |Max(R)|+1$. Also, $|K|=|Max(R)|$ iff $R=\mathbb Z_2$ or $R\cong\mathbb Z_2\times \mathbb Z_2$.  }

\vs{3mm}\Proof (1) Clear.

(2). Assume to the contrary that there exists some nonzero $u\in \mf m\cap\mf n\cap K$. Since $\mf m+\mf n=R$, we have $x\in \mf m,y\in\mf n$ such that $x+y=1$. Clearly,  $x\not=y$ and hence we can assume $x\in K$. Then it follows by (1) that $x=u\in \mf n$, a contradiction.

(3) If $R$ is a local ring, then either $R\cong \mathbb Z_2$ or $R$ has at least two units by assumption. The result holds in this case. In the following, we assume that $R$ is non-local. If there is a maximal ideal $\mf m$ with
$\mf m\seq D$, then by the proof of \cite[Theorem 2.1]{JM}, $R\cong \mathbb Z_2\times \mathbb F$ for some field $\mathbb F$. In this case, $|K|\ge |Max(R)|$ and, $|K|=|Max(R)|$ iff $\mathbb F=\mathbb Z_2$. If no maximal ideal is contained in $K$, then  each maximal ideal has exactly one vertex in $K$ by $(1)$. By (2), we have $|K|\ge |Max(R)|+1$.\qed

\vs{3mm}We remark that Lemma 2.1 (3) is the best possible result. For instance, $|K|=|Max(R)|+1$ holds for the local ring $\mathbb Z_4$ and $\prod_{1}^3\mathbb Z_2$. For $R=\prod_{1}^3\mathbb Z_2$,  we draw its co-maximal graph $\Om(R)$ in Figure 1, in which $\Om(R)=K_1+K_1+H$.

\vs{3mm}The following is a result of \cite{JM}:

\vs{3mm}\nin{\bf Lemma 2.2.} (\cite[Theorem 2.1]{JM}) {\it For a commutative ring $R$ which is non-local, if the co-maximal graph $\Om(R)$ is a split graph and $R\not\cong \mathbb Z_2\times \mathbb F$ for any field $\mathbb F$, then the characteristic of $R$ is two, $R$ has exactly three maximal ideal $\mf m_i$, and  $K\cap \mf m_i=\{x_i\} $ for all $i$, where each $x_i$ is idempotent and for every invertible element $u$ of $ R$, $ux_i=x_i$.   }

\vs{3mm} We now prove give an alternative proof to the main result of \cite{JM}:

\vs{3mm}\nin{\bf Theorem 2.3.} (\cite[Theorem 2.3 ]{JM}) {\it For a commutative ring $R$, the co-maximal graph $\Om(R)$ is a split graph if and only if $R$ is one of the following: a local ring, $\mathbb Z_2\times \mathbb Z_2\times \mathbb Z_2$, $\mathbb Z_2\times \mathbb F$ for some field $\mathbb F$.}

\Proof We only need  prove the necessary part. Assume that $R$ is not a local ring and $R\not\cong \mathbb Z_2\times \mathbb F$ for any field $\mathbb F$. Let $G=\Om(R)$. Assume further that the co-maximal graph $G$ is a split graph with $V(G)=K\cup D$. By Lemma 2.2, $R$ has exactly three maximal ideals $\mf m_i$ ($1\le i\le 3$) such that $\mf m_i\cap K=\{x_i\}$, where $x_i^2=x_i$ and $ux_i=x_i$ holds for every $u\in U(R)$. Also, $R$ has characteristic 2. For any $1\le i\not=j\le 3$, assume $1=rx_i+sx_j$ with $rx_i\in K$. By Lemma 2.1 we have $x_i=rx_i$ and hence $1-x_i=sx_j$. It follows that $(1-x_i)(1-x_j)=0$.
Notice that $1+x_1x_2x_3\in U(R)$ since $x_1x_2x_3\in J(R)$, it follows that $x_1x_2x_3=0$. Then  $(1-x_i)(1-x_j)=0$ implies $x_k+x_ix_k+x_jx_k=0$ for any re-arrangement $i,j,k$ of $1,2,3$.  Therefore $x_1+x_2+x_3=2\cdot f(x_1,x_2,x_3)=0$. This shows that $1-x_1,1-x_2,1-x_3$ is a complete set of orthogonal idempotent elements of $R$.

Let $R_i=R(1-x_i)$. Then $R=\sum_{i=1}^3R_i\cong R_1\times R_2\times R_3$, where each $R_i$ is a local ring with a unique maximal ideal $\mf n_i$, since $R$ has exactly three maximal ideals. For any $r(1-x_1)\in \mf n_1$, we have $$1+r(1-x_1)=[(1-x_1)+r(1-x_1)]+(1-x_2)+(1-x_3)\in U(R).$$ Since $ux_2=x_2$ for all $u\in U(R)$, it follows that $r(1-x_1)\cdot x_2=0$. Now apply the fact $(1-x_1)(1-x_2)=0$, one derives $r(1-x_1)=0$. Hence $\mf n_1=0$ and each $R_i$ is a field. Thus $R$ is a direct product of three fields.

Since $x_1=x_2+x_3$, it follows that $$\mf m_1=R_2+R_3,\,\mf m_2=R_1+R_3,\,\mf m_3=R_1+R_2.$$
Then in the decomposition $R=K\cup D$, we deduce from Lemma 2.1 (1) that $$K=(U(R_1)+U(R_2)+U(R_3))\cup \{x_1,x_2,x_3\},\,D=R\setminus K.$$
Now we claim that each $R_i$ is isomorphic to $\mathbb Z_2$ and hence,  $R\cong \mathbb Z_2\times\mathbb Z_2\times\mathbb Z_2 .$ In fact,  assume $R=R_1\times R_2\times R_3$. Then $$K=(U(R_1)\times U(R_2)\times U(R_3))\cup \{(0,1,1),(1,0,1),(1,1,0)\}.$$ If $|R_1|>2$, then there exists a nonzero element $1\not=v_1\in R_1$. Then both $z=(v_1, 0,1)$ and $e=(0,1,0)$ are in the independent subset  $D$, contradicting $Rz+Re=R$. This completes the proof.\qed

\vs{3mm}Recall a convenient construction from graph theory, {\it the sequential sum} $$G_1+G_2+\cdots+G_r$$ of a sequence of graphs $G_1,G_2,\ldots,G_r$. We illustrate the construction in Figure 1 for the sequence of graphs $K_1, K_1, H$, where $H$ is a triangle $K_3$ together with three end vertices adjacent to distinct vertices.

\setlength{\unitlength}{0.10cm}
\begin{figure}[h]
 \centering
\input{1.LP}
\end{figure}

\vs{3mm}\nin{\bf Corollary 2.4. } {\it A finite split graph is a co-maximal graph of a non-local commutative ring if and only if $G$ is one of the following: (1) The sequential sum $K_1+K_{p^n-1}+K_{1,p^n-1}$ for some prime number $p$. (2) The sequential sum $K_1+K_1+H$ in Figure 1, where $H$ is a triangle $K_3$ together with three end vertices adjacent to distinct vertices.}

\vs{3mm}The co-maximal graph of a field is certainly a complete graph. For a finite local ring $(R,\mf m)$ which is not a field, assume $|\mf m|=p^m$.  Then the co-maximal graph of $R$ is $K_{p^n-p^m}+D_{p^m}$, where $p^n\ge2 p^m$ and $D_{p^m}$ is a discrete graph with $p^m$ vertices.

\vs{4mm}\nin{\bf 3. The subgraph  $\G(R)$}

\vs{4mm}\nin As noticed by Maimani et.al. in \cite{MSSY}, the main part of the co-maximal graph $\Om(R)$ is the subgraph $\G(R)$ induced on the vertex subset $R\setminus (U(R)\cup J(R))$. In fact, we have the following facts:

\vs{3mm}(O1) A vertex in $U(R)$ is adjacent to every vertex of $\Om(R)$,  while an element of $J(R)$ only connects to units of $R$. In fact, there is a sequential sum decomposition $$\Om(R)=J(R)+U(R)+\G(R).$$

(O2) $\G(R)$ is empty if and only if $R$ is a local ring, i.e., a commutative ring with a unique maximal ideal.

\vs{3mm}\cite{MSSY,wang} studied this subgraph and obtained many interesting results. In particular, it is proved that the graph is connected with diameter less than or equal to three (\cite[Theorem 3.1]{MSSY}), that the girth of the graph is less than or equal to four (\cite[Corollary 3.8]{wang}). We include a detailed proof for the following fundamental property of graphs related to algebras:

\vs{3mm}\nin{\bf Theorem 3.1.(\cite[Theorem 3.1]{MSSY} )}  {\it
The graph $\G(R)$ is connected with diameter less than or equal to three.}

\Proof For any $a\in R$, set $S_a=\{\mf m\in Max(R)\,|\,a\in \mf m\}$. Then for each $ a\not\in J(R)$, $S_a\subset Max(R)$. For distinct $a,b\in R$, we claim

(1) $ab\in J(R)$ iff $Max(R)=S_a\cup S_b$,

(2) $Ra+Rb=R$ iff $S_a\cap S_b=\eset.$

Now for distinct  $a,b\in R\setminus (U(R)\cup J(R))$,  if $ab\not\in J(R)$, then there exists $x\in R\setminus (U(R)\cup J(R))$ such that $Rab+Rx=R$. Then clearly there is a path $a-x-b$ in $\G(R)$ and hence $d(a,b)\le 2$. If $ab\in J(R)$, then take any $y\in R\setminus (U(R)\cup J(R))$ such that $Ra+Ry=R$. We claim that $by\not\in J(R)$ and it will follow that $d(a,b)\le d(a,y)+d(y,b)\le 3.$ In fact, assume to the contrary that $by\in J(R)$. Then we have
$$S_b\cup (S_y\setminus S_b)=S_b\cup S_y=Max(R)=S_b\cup (S_a\setminus S_b).$$ It follows that $(S_y\setminus S_b)=S_a\setminus S_b=\eset$ since $S_a\cap S_y=\eset$. Then $b\in J(R)$, a contradiction. This completes the proof.\qed

\vs{3mm}By \cite[Theorem 3.9(1)]{wang}, the clique number $\om(\G(R))$ is infinite whenever $J(R)=0$ and the ring $R$ is indecomposable. It could be used to sharpen \cite[Theorem 2.2]{MSSY}, as the following theorem shows.

\vs{3mm}\nin{\bf Theorem 3.2. } {\it For any non-local ring $R$, let  $G=\G(R)$. Then the following are equivalent:}

 (1) {\it $G$ is a  bipartite graph.}

 (2) {\it $G$ is a complete bipartite graph.}

 (3) {\it $R$ has exactly two maximal ideals.}

 (4) {\it $R/J(R)\cong \mathbb K_1\times \mathbb K_2$, where each $\mathbb K_i$ is a field.}

\Proof $(3)\Lra (2)$: Assume that $\mf m_1,\mf m_2$ are the maximal ideals of $R$. Then $\G(R)$ is a complete bipartite graph with two partitions $\mf m_1\setminus \mf m_2$ and $\mf m_2\setminus \mf m_1$.

$(2)\Lra (3)$: Assume that $G$ is a complete bipartite graph with vertex partition $V(G)=V_1\cup V_2$. Then for any maximal ideal $\mf m$, $\mf m\setminus J(R)$ is entirely contained in a single partition.
Assume $\mf m_i\setminus J(R)\seq V_i$. If $R$ has a third maximal ideal $\mf n$, then $\mf n\setminus J(R)\seq V_i$ for some $i$. This is impossible since $\mf n+\mf m_i=R$.

$(2)\Lra (1)$ and $(4)\Lra (3)$: Clear.

$(1)\Lra (3)$: Clearly, there is no loss to assume $J(R)=0$.

If $R$ has only a finite number of maximal ideals, say $\mf m_i,\, 1\le i\le r$, set $\mf n_i=\mf m_i\setminus (\cup_{j\not=i}\mf m_j)$.
Then $\mf n_i\not=\eset$, and each vertex in $\mf n_i$ is adjacent to all vertex in $\mf n_j$. This implies $\omega(G)\ge |Max(R)|$ and hence $\omega(G)= |Max(R)|$, when $R$ has only a finite number of maximal ideals. If further $G$ is a bipartite graph, then $r=2$, i.e., $R$ has exactly two maximal ideals.

In the following, assume that $R$ has infinitely many maximal ideals, and we proceed to prove $\omega(G)=\infty$. Assume that $R$ has a non-trivial idempotent $e$. Then $R=eR\times (1-e)R$. If $eR$ has no nontrivial idempotent element, then $\omega(\G(eR))=\infty$ and hence $\omega(G)=\infty$. Then assume that both $e$ and $1-e$ are non-primitive idempotents. By induction, for any integer $s\ge 1$, there exist non-trivial corner rings $R_{s,i}$ of $R$ such that $R=\prod_{j=1}^sR_{s,j}$. Set $$f_1=(0,1,\cdots, 1),\cdots, f_2=(1,0,1,\cdots, 1),\cdots, f_s=(1,\cdots, 1,0).$$ Clearly $\{f_j\}_{1\le j\le s}$ is a clique in $\G(R)$, and thus $\omega(G)=\infty$.

$(3)\Lra (4)$: This follows from the Chinese Remainder Theorem. \qed

\vs{3mm} Notice that the proof together with \cite[Theorem 3.9(1)]{wang} actually gives an alternative proof to the fact that $\omega(\G(R))=|Max(R)|$ whenever $R$ is not a local ring.

\vs{3mm}Recall that a ring $R$ is called an {\it exchange ring} if the left module $_RR$ has the exchange property, see \cite{WuXu} and the included references for details. Recall that idempotents can be lifted modulo every ideal of an exchange ring $R$. The class of exchange rings include artinian rings, semiperfect rings and clean rings, the rings in which each element is a sum of an idempotent and a unit. Recall that for a ring $R$ with all idempotents central in $R$, $R$ is clean iff $R$ is an exchange ring. For commutative clean rings $R$, we have

\vs{3mm}\nin{\bf Corollary  3.3. } {\it For any  commutative non-local exchange ring $R$, let  $G=\G(R)$. Then
$G$ is a  bipartite graph iff  $G$ is a complete bipartite graph, if and only if $R\cong  R_1\times R_2$, where each $ R_i$ is a local ring.}

\vs{3mm} A graph $G$ is called {\it totally disconnected} if the edge set $E(G)$ is empty. By Theorem 3.1, we have the following observation and hence Theorem 3.4:

\vs{3mm}(O3) $\G(R)$  is totally disconnected iff it is an empty graph, if and only if $R$  is a local ring.

\vs{3mm}\nin{\bf Theorem 3.4.(\cite{MSSY,wang} )}  {\it  For a ring $R$, let $G=\G(R)$. Then the following are equivalent:}

(1) {\it $G$ is a refinement of a star graph, i.e., $G$ has at least two vertices, and there exists a vertex in $G$ which is adjacent to every other vertex.}

(2) {\it $G$ is a tree, i.e., $G$ is nonempty, connected and contains no cycles.}

(3) {\it $G$ is a star graph.}

(4) {\it $R$ is isomorphic to  $\mathbb Z_2\times \mathbb F$ for some field $\mathbb F$.}

\Proof $(1)\Llra (3)\Llra (4)$:  This is contained in \cite[Corollary 2.4(2)]{MSSY}.

$(2)\Llra (4)$: This is contained in \cite[Theorem 3.5, Corollary 3.6]{wang}. \qed

\vs{3mm}\nin{\bf Corollary 3.5. } {\it For a finite simple graph $G$ with $|G|\ge 2$, assume $G=\G(R)$ for some ring $R$. Then the following are equivalent:}

(1) {\it $G$ is a refinement of a star graph.}

(2) {\it $G$ is a tree, i.e., $G$ contains no cycle.}

(5) {\it $G=K_{1,p^n-1}$ for some prime number $p$ and some positive integer $n$.}

\vs{3mm}By the result of section 2 and the results of \cite{MSSY,wang}, it is natural to ask the following question: For what rings $R$ is $\G(R)$ a split graph? Notice that for distinct
maximal ideals $\mf{ m,n}$ (if exist) and $x\in \mf m,y\in \mf n$, $Rx+Ry=R$ implies $x,y\in R\setminus (U(R)\cup J(R))$. Then a careful check to the proofs of Lemma 2.1, Lemma 2.2 and Theorem 2.3 shows the following:

\vs{3mm}\nin{\bf Theorem 3.6. } {\it For any ring $R$, the following statements are equivalent:}

(1) {\it The co-maximal graph $\Om(R)$ of $R$ is a split graph.}

(2) {\it The subgraph $\G(R)$ is either empty or a split graph.}

(3) {\it $R$ is one of the following: a local ring, $\mathbb Z_2\times \mathbb Z_2\times \mathbb Z_2$, $\mathbb Z_2\times \mathbb F$ for some field $\mathbb F$. }

\vs{3mm} Notice that a split graph $G$ is isomorphic to $\G(R)$ for some finite ring $R$ iff either $G$ is a star graph $K_{1,p^n-1}$ for some prime number $p$, or $G$ is the triangle together with three end vertices adjacent to distinct vertices (see Figure 1).

\vs{3mm} The works of \cite{MSSY,wang} show that the graph $\G(R)$ has many properties which the zero-divisor graph of a ring (or a semigroup) already have. Recall from \cite{DMS,Mulay} that {\it the core of a zero-divisor graph $G$ is always a union of triangles and squares,  and a vertex in $G$ is either an end vertex or a vertex of the core}. Recall that  the core of a graph $G$ is by definition the subgraph induced on all vertices of cycles of $G$. In the final part of this section, we will show that the graph $\G(R)$ has the same property.

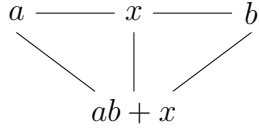
\vs{3mm}\nin{\bf Lemma 3.7.} {\it For any path $a-x-b$ in the graph $\G(R)$, if $ab+x$ is not a unit of $R$, then the path is contained  in a subgraph isomorphic to $K_1+K_2+K_1$. }

\Proof  The given condition implies $Rx+Rab=R$. Hence $ab\not\in J(R)$ and in particular, $ab+x\not=x$.  Furthermore, we have $$R(ab+x)+Ra=R=R(ab+x)+Rb.$$ Now assume $ab+x\not\in U(R)$. Then $ab+x\not\in U(R)\cup J(R)$. If $ab+x\not\in \{a,b\}$, then there is a subgraph $K_1+K_2+K_1$ in $\G(R)$ which contains the path $a-x-b$, see Figure 2. Since $R(ab+x)+Ra=R=R(ab+x)+Rb$, it follows that $ab+x\not\in \{a,b\}$, and this completes the proof. \qed

\setlength{\unitlength}{0.13cm}
\begin{figure}[h]
 \centering
\input{2.LP}
\end{figure}

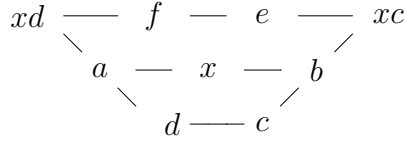
\vs{3mm}\nin{\bf Lemma 3.8.} {\it If a vertex $x$ of $\G(R)$  is in a cycle of five vertices,  then $x$ is  in either a triangle or  a rectangle.}

\Proof Assume $a-x-b-c-d-a$ is a cycle in $\G(R)$, and $x$ is not in any triangle. We proceed to verify that $x$ is in a square.

\setlength{\unitlength}{0.12cm}
\begin{figure}[h]
 \centering
\input{4.LP}
\end{figure}

In fact, if $xv\not=x$ for some $v\in U(R)$, then there is a square $a-x-b-vx-a$ in $\G(R)$. Thus in the following, we assume $uy=y,\,\forall u\in U(R)$ where $y\in \{a,x,b\}$. Then by Lemma 3.7, we can assume further $$ab+x=1,\,xc+b=1,\,xd+a=1\quad (*).$$ Notice that $xd$ is adjacent to $a$, and $xc$ is adjacent to $b$ in $\G(R)$. If $xd=x$,
then by $(*)$ we have $a=ab$. Then there is a rectangle $a-x-b-d-a$ in $\G(R)$. Therefore assume $xd\not\in \{x,d\}, \, xc\not\in \{x,c\}. $

Since $d(xc,xd)\le 3$, we can assume there is a path $xc-e-f-xd$. If $e\not=b$, then there is a rectangle $x-e-xc-b-x$. If $e=b$, then $f\not\in \{x,a\}$
 and hence by $(*)$ there is a cycle $a-x-b-f-xd-a$. Then there is a rectangle $a-x-f-xd-a$ and a triangle $b-x-f-b$. This completes the verification. \qed


\vs{3mm} \vs{3mm}\nin{\bf Lemma 3.9.} {\it Let $G=\G(R)$ and assume that $G$ contains a cycle. Then for any path $a-x-b$ in the core of $G$, $x$ is  in a cycle $C_n$ with  $3\le n\le 5$. }

\Proof Assume $$\cdots -f-a-x-b-g-\cdots $$ is a cycle in $\G(R)$. We proceed to verify that  $x$ is in a cycle $C_n$ with $3 \le n\le 5$.  Since $Rab+Rfb=Rb$, we have $$Rb=Rab+(Ra+Rx)fb=Rab+Rxfb.$$
Then $Rb+R(x-xfb)=Rab+Rx=R.$ It follows that $x-xfb\not\in U(R)\cup J(R).$  Consider $xfb.$ If  $xfb=0$, then it follows that $Rb=Rab+Rxfb=Rab.$ In this case, there is a rectangle $a-x-b-g-a$. In the following, assume  $xfb\not=0$, i.e., $x-xfb\not=x$, and let $y=x-xfb.$ Without loss of generality, assume $d(y,f)=3$ and  $y-d-e-f$ is a path from $y$ to $f$. If further $d\not=b$, then there is a rectangle $b-x-d-y-b$.  If  $d=b$, then there is a cycle $C_5:$ $f-a-x-b-e-f$. This completes the proof \qed

\vs{3mm} \vs{3mm}\nin{\bf Theorem 3.10.} {For a ring $R$, let $G=\G(R)$ and assume that $G$ contains a cycle. Then the core of $G$ is a union of triangles and rectangles, and every vertex of $G$ is either an end vertex or a vertex of the core.}

\Proof The first statement follows from Lemmas 3.9 and 3.8. The second statement follows after an argument similar to the proof of \cite[Theorem 1.5]{DMS}. We omit the details here.\qed

\vs{4mm}\nin{\bf 4. A retract $\G_r(R)$ of $\G(R)$}

\vs{3mm}\nin For simple graphs $G$ and $H$, recall that  $G\os{\vi}{\to} H$ if there exists a map $\varphi: G\to H$  such that for distinct $u,v\in V(G)$, $u-v$ in $G$ implies $\varphi(u)\not=\vi(v)$ and $\vi(u)-\vi(v)$ in H. Such a map is called a {\it graph homomorphism}. If $H$ is a subgraph of $G$, $G\os{\vi}{\to} H$ and the restriction of $\vi$ on $H$ is an identity, then $H$ is called a {\it retract } of $G$, see Figure 4 for an example. If $G$ has no proper retract, then $G$ is called a {\it core graph} (e.g., $K_n$ is a core graph.). Notice that the core of a graph needs not be a core graph. Recall that the {\it chromatic number} $\chi(G)$ of $G$ is the least positive integer $r$ such that $G\to K_r$. It is the least number of colors needed for coloring the vertices of $G$ in such a way that no two adjacent vertices have a same color. The girth of
a graph $G$, denoted by $g(G)$, is the length of the minimal cycle in $G$.

In this section,
we introduce a new graph which is a retract of $\G(R)$ and we study this new graph.

\vs{3mm}\nin{\bf Definition 4.1.} {\it For a ring $R$ and any $x\in R$, let $\ol{x}=Rx$. Construct a simple graph in the following and denote it as $\G_r(R)=G$: }
$$V(G)=\{\ol{x}\,|\,x\in R\setminus (U(R)\cup J(R))\},$$
$$E(G)=\{\{\ol{x},\ol{y}\}\,|\, \text{$\ol{x}\not=\ol{y}\in V(G)$, and $Rx+Ry=R$} \}$$

Clearly, this graph has less vertices and less edges than that of the graph $\G(R)$ in general, see Figure 4. More precisely,

\vs{3mm}\nin{\bf Proposition 4.2.} (1) {\it For a ring $R$, the graph $\G_r(R)$ is a retract of $\G(R)$.  }

(2) {\it For any ideal $I$ of $R$ contained in $J(R)$, $\G(R/I)$ is a retract of $\G(R)$.}

{\it Furthermore},

(3)  {\it If the girth of $\G(R)$ is three, then so is the girth of $\G_r(R)$.}

(4) {\it $\G_r(R)$ is connected with diameter less than or equal to three.}

(5) {\it $\omega(\G_r(R))=\omega(\G(R))$ and  $\chi(\G_r(R))=\chi(\G(R)).$}

(6) {$\G_r(R)$ and $\G(R)$ are homomorphically equivalent to a same unique core graph.}

\Proof  $(1)$ Clearly, $\vi: x\mapsto \ol{x}$ is a surjective graph homomorphism. For any $\ol{x}\in V(\G_r(R))$, fix a vertex $y_x$ in $Rx$ to obtain a subgraph of $\G(R)$. Then the graph $\G_r(R)$ is isomorphic to the subgraph of $\G(R)$. Thus $\G_r(R)$ is a retract of the graph $\G(R)$. In a similar way, one checks $(2)$. All the remaining results follow from $(1)$. We check $(3)$ and $(5)$ in the following.

$(3)$ The result is clear by the definition of $\G_r(R)$. In fact, if $G\to H$ and $G$ has odd girth, then it is known that $g(G)=g(H)$, since any cycle with odd length is a core graph.

(5) Since a composition of graph homomorphisms is still a graph homomorphism, clearly $\chi(\G_r(R))\ge \chi(\G(R))$. On the other hand, since $\G_r(R)$ is isomorphic to a subgraph of $\G(R)$, for any graph homomorphism $\psi: \G(R)\to K_s$, the restriction  $\psi |: \G_r(R)\to K_s$ is also a graph homomorphism. This shows the second assertion of (5). The first equality holds by the definition of a retract of a graph.
\qed

\vs{3mm}
By Theorem 3.2 and Proposition 4.2, we have

\vs{3mm}\nin{\bf Corollary 4.3. } (1) {\it For any non-local ring $R$, $\chi(\G_r(R))\ge |Max(R)|$.}

(2) {\it $\G_r(R)$ is a complete bipartite graph iff $R$ has exactly two maximal ideals, iff $\G_r(R)$ is a bipartite graph.}

\vs{3mm} By \cite[Corollary 3.8(1)]{wang}, if $H=\G(R)$ contains a cycle, then the girth $g(H)\le 4.$ Compared with Proposition 4.2(3), we have the following example for the $g(H)=4$ case:

\setlength{\unitlength}{0.12cm}
\begin{figure}[h]
 \centering
\input{3.LP}
\end{figure}

\vs{3mm}\nin{\bf Example 4.4.} Consider the ring $R=\mathbb Z_{12}$. For this ring, $g(\G(R))=4$, while $g(\G_r(R))=\infty$. We draw the two graphs in Figure 4.

\vs{3mm}\nin{\bf Theorem 4.5. } {\it For any commutative ring $R$ which is not a local ring, let $G=\G(R)$. Then the following numbers are identical:}

(1) {\it The chromatic number $\chi(G)$;}

(2) {\it The clique number $\om(G)$;}

(3) {\it The cardinal number of $Max(R)$;}

(4) {\it $\chi(\G_r(R))$; }

(5) {\it $\om(\G_r(R))$.}

\Proof (1) First, recall from \cite[Theorem 2.3]{SB} that $$\chi(\Om(R))=|Max(R)|+|U(R)|.$$ Then consider the following decomposition of the co-maximal graph into sequential sums of three subgraphs
 $$\Om(R)=J(R)+U(R)+\G(R),$$
 where $J(R)$ is a discrete subgraph while $U(R)$ is a complete subgraph. Since $|J(R)|\le |U(R)|$, we have $$\chi(\Om(R))=\chi(\G(R))+|U(R)|.$$ By \cite[Theorem 3.9(2)]{wang}, $\om(\G(R))=|Max(R)|$. Then $$\chi(\G(R))=|Max(R)|=\om(\G(R))$$

 (2) By Lemma 4.2, we have $\chi(G)=\chi(\G_r(R))\ge \om(\G_r(R))=\om(G)$, and the result follows from $(1)$. \qed

\vs{3mm}Now we study the interplay between the graph structure of $\G_r(R)$ and the algebraic property of $R$. First we have

\vs{3mm}\nin{\bf Proposition 4.6.} {\it For a ring $R$, let $G=\G_r(R)$. Then the following are equivalent:}

(1) {\it $G$ is a refinement of a star graph.}

(2) {\it $G$ is a star graph.}

(3) {\it $R\cong \mathbb F\times T$, where $\mathbb F$ is a field and $T$ is a local ring. }

\Proof We  only need prove $(1)\Lra (3)$.
Assume $\ol{x}$ is adjacent to every other vertex in $\G_r(R)$. Then $\ol{x}=\ol{x^2}$. Assume $x=rx^2$. Then $rx=(rx)^2$ and $\ol{x}=\ol{rx}$. So we can assume at the start that
$x$ is a nontrivial idempotent. Let  $Rx=T,\,\mathbb F=R(1-x)$. Then $T$ is a maximal ideal of $R$ and $R=\mathbb F\times T$. Thus $\mathbb F$ is a field. Now $J(R)=J(T)$, and for any $y\in T\setminus J(T)$, we have $T=Ty$. Thus $T$ is a local ring with the maximal
ideal $J(R)$. \qed

\vs{3mm} Notice that for a finite local ring $(T,\mf m)$, if $\mf m=Tx_1\cup \ldots \cup Tx_r$ in which $Tx_i\not=Tx_j$, then $\G_r(R)=K_{1,r}$.

\vs{3mm}\nin{\bf Corollary 4.7. } {\it For any commutative ring $R$, $diam(\G_r(R))=1$ iff $R\cong \mathbb F_1\times \mathbb F_2$, where each $\mathbb F_i$ is a field. }

\Proof If $R\cong \mathbb F_1\times \mathbb F_2$, then clearly $diam(\G_r(R))=1$. Conversely, assume $diam(\G_r(R))=1$. Then $\G_r(R)$ is a complete graph. By Proposition 4.6, we have $R\cong \mathbb F \times T$, where $\mathbb F$ is a field and $T$ is a local ring. If $T$ is not a field, then take any nonzero $x\in J(T)$. Then $(1,x)\in R\setminus U(R)\setminus J(R)$, and $\ol{(1,0)}$ is not adjacent to $\ol{(1,x)}$ in $\G_r(R)$, a contradiction. This completes
the proof.\qed

\vs{3mm}The following result follows from Corollary 4.3, Proposition 4.6 and the proof of \cite[Lemma 3.2, Propositiom 3.3]{MSSY}:

\vs{3mm}\nin{\bf Proposition 4.8. } {\it For any commutative non-local ring $R$, $diam(\G_r(R))=2$ iff  one of the following conditions holds:  }

(1) {\it $J(R)$ is a prime ideal of $R$.}

(2) {\it $R$ has exactly two maximal ideals, and $R\not\cong \mathbb F_1\times \mathbb F_2$ for any fields $\mathbb F_i$. }

\vs{3mm}\nin{\bf Corollary 4.9. } {\it For any commutative non-local ring $R$, the graphs $\G(R)$ and $\G_r(R) $ have a same diameter iff  $R\not\cong \mathbb F_1\times \mathbb F_2$ for any fields $\mathbb F_i$ but $\mathbb F_1=\mathbb F_2=\mathbb Z_2$.}

\vs{3mm} By \cite[Theorem 3.9(2)]{wang}, Theorem 4.5  and Corollary 4.7, if $diam(\G_r(R))=2$ (respectively, $diam(\G(R))=2$), then either $\G_r(R)$ (respectively, $\G(R)$) is a complete bipartite graph or its clique number is infinite.

\vs{3mm}At the end of the paper, we pose the following problem:

\vs{3mm}\nin{\bf Question 4.10. } {\it Which rings $R$ have the property that $\G(R)$ is a generalized split graph? Which  rings $R$ have the property that $\G_r(R)$ is a (generalized) split graph?}

\vs{3mm}\nin Recall from \cite{LW} that a simple graph $G$ is a {\it generalized split graph} if
$$V(G)=K\cup D, K\cap D=\eset,$$ where the induced subgraph on $K$ (resp., on $D$) is a core graph (respectively, a discrete graph). Notice that $\Gamma(\mathbb Z_{12})$ is not a generalized split graph, while $\G_r(\mathbb Z_{12})$ is a split graph.

\end{document}